\documentclass[amstex,12pt,amssymb]{article}
\usepackage{mathtext}
\usepackage{amsmath}
\usepackage{amssymb}
\usepackage{amsxtra}
\usepackage{latexsym}
\usepackage{ifthen}

\usepackage{indentfirst}

\usepackage{color}

\def\DD{P\ r\ o\ o\ f}

\makeatletter
\def\thmstyle{\it} 
\def\@begintheorem#1#2{\it \trivlist \item[\hskip
        \labelsep{\bf #1\ #2.}]\thmstyle}
\def\@opargbegintheorem#1#2#3{\it \trivlist \item[\hskip
        \labelsep{\bf #1\ #2\ (#3).}]\thmstyle}
\makeatother

\begin{document}

\renewcommand{\baselinestretch}{1.25}
\newcommand {\beq}{\begin{equation}}
\newcommand {\eeq}{\end{equation}}

\newtheorem{theorem}{Theorem}[section]
\newtheorem{lemma}[theorem]{Lemma}
\newtheorem{proposition}[theorem]{Proposition}
\newtheorem{corollary}[theorem]{Corollary}
\newtheorem{conjecture}[theorem]{Conjecture}
\newtheorem{definition}[theorem]{Definition}
\newtheorem{remark}[theorem]{Remark}
\newtheorem{claim}[theorem]{Claim}
\newtheorem{cons}[theorem]{Consequence}

\title{Limit theorems for chains with unbounded variable length memory
which satisfy Cramer condition
\footnote{The research is supported by RSF according to the research project 18-11-00129, AL and AY
thank FAPESP grant 2017/20482-0, AY also thanks CNPq and FAPESP grants 301050/2016-3 and 2017/10555-0,
respecively. It is part of USP project Mathematics, computation, language and the brain, FAPESP project
Research, Innovation and Dissemination Center for Neuromathematics
grant 2013/07699-0}
 \author{\bf A. Logachev, A. Mogulskii,  A. Yambartsev
 }
}\maketitle

\begin{abstract} Here we obtain the exact asymptotics for large and 
moderate deviations, strong law of large numbers and central limit theorem for
chains with unbounded variable length memory.
\end{abstract}


\section{Instroduction}

Let $\mathcal{A}=\{0,1, \dots, d\}$ be a set of $d$ symbols (characters) -- an alphabet. Here we
consider a class of chains $\mathbf{r} = (r_i, i\in \mathbb{Z}) \in \mathcal{A}^{\mathbb{Z}}$
which is a special case of so-called \textit{chains with unbounded
variable length memory}. These chains began to be actively studied after they were first introduced
by Jorma Rissanenom \cite{R} as an economic and universal way of data compression. 
The short and simple introduction to these processes can be found, for example, in  \cite{GL2008}.
They are used for modeling data in computer science 
\cite{R}, in biology \cite{BY2001},
\cite{L2006},  in neuro-biology  \cite{DO2016}, \cite{DGLO2018}, in linguistics \cite{4GL}. 
We do not pretend for at least some sort of full review on these chains, and we restrict ourselves to some 
papers known to us by the activity of the research group \textit{NeuroMat} under the guidance of Prof. A.Galves.

If we interpret $\mathbb{Z}$ as a discrete time, then one can imagine such chains as a successive 
in time attribution of a character from the alphabet $\mathcal{A}$ with probability which depends on the
past (existing sequence of characters), or, more precisely, it depends on a part of the past, a context.
As a consequence such dependence can be represented as a context tree, where each vertex represent a 
context, and each vertex is associated with a probability distribution (on $\mathcal{A}$) of a new 
character. Markov chain with the state space $\mathcal{A}$ is a particular case of these chains. In this
case the context tree has height 1 -- we should know only the last character in order to decide the 
distribution of a new character. 
The existence of a stationary measure on $\mathcal{A}^{\mathbb{Z}}$ compatible with a family of 
transition probabilities determined by a context tree is the question which naturally arises. These
question was answered and the short review can be find \cite{GL2008}, where also methods of statistical
inferences of context tree were provided.

In \cite{Gallo2011} the perfect simulation for such processes was constructed. The success of the algorithm
(if the perfect simulation stops in a finite time) depends directly from the existence of a (finite) renewal
time moment -- the moment when the next and successive attributions of characters do not depend on the past.
In the same paper the connection between renewal processes and the chains with variable length memory was
established. This indicates to us that the large deviations results for this sort of chain can be obtained
using the regeneration structure and the recently published results on large deviations for the 
renewal processes \cite{Mog1}.



Despite the fact that for a complete definition and description of chains with unbounded variable 
length memory we need to introduce the notion of a context tree, in this article we restrict ourselves 
to an alternative description of the chain, because from the very beginning we will consider a particular
case of such chains. Let us fix an initial configuration

$$
r(0)=\{r_i(0)\}_{-\infty< i \leq 0},
$$
where $r_i(0)$, $-\infty< i \leq 0$ take values from the alphabet $\mathcal{A}$, which is binary
alphabet $\mathcal{A}:=\{0,1\}$ in our case. We set the configuration change rule. Remind that
at any time step we write exactly one character from the alphabet $\mathcal{A}$ at the end of 
existing sequence without changing it.

In order to set the transition rules, let first fix a number $v\in \mathbb{N}$ (one of the parameters of 
the chain). Consider the set of all words from $v$ characters ending with the character $1$; 
the total number of such words is $2^{v-1}$; for any such word we assing their own order number $j$,
$1\leq j\leq 2^{v-1}$. Let us set the sequence of positive numbers on this set of words
$$
\mathbf{p}_{kj} \in (0,1), \ \ \ k\in 0\cup\mathbb{N}, \ \ \ 1\leq j\leq 2^{v-1}.
$$

Now we are ready to describe the rules (transition probabilities) of adding a character from the right:
let we have a configuration on the $n$-th step
$$
r(n)=\{r_i(n)\}_{-\infty< i \leq n}.
$$
Denote
$$
m_n:=\sup\{-\infty< i \leq n:r_i(n)=1\}.
$$
%
Then, at the next step the configuration $r(n)=\{r_i(n)\}_{-\infty< i \leq n}$ jump to the 
configuration
$$r(n+1)=\{r_i(n+1)\}_{-\infty< i \leq n+1}$$
by writing from the right a character $1$, $r_{n+1}(n+1)=1$, with probability 
$\mathbf{p}_{k_{n}j_n},$
where $k_{n}=n-m_n,$ $j_n$ -- the number of word, which forms the sequence 
$$
r_{m_n-v+1}(n),r_{m_n-v+2}(n), \dots, r_{m_n}(n).
$$
Thus, the character $0$, $r_{n+1}(n+1)=0$, will be adding with probability $1-\mathbf{p}_{k_{n}j_n}$.
Note that the previous sequence do not change:
$$
\{r_i(n+1)\}_{-\infty< i \leq n} \equiv r(n).
$$
Thus, the probability of the attributed character $r_{n+1}(n+1)$ depends on
\begin{enumerate}
 \item[1)] the distance to the nearest character $1$;
 \item[2)] the word from $v-1$ letters, which stand on the left from this $1$.
\end{enumerate}

It is obvious that the random sequence $\{r_{n}(n)\}$ is not Markov chain, because
the transition probability from $r_{n}(n)$ to $r_{n+1}(n+1)$ can depend, generally speaking,
not only on the character $r_{n}(n)$, but also on values 
$r_{n-j}(n-j)$ for arbitrarily large $j\geq 0$.

For such defined process $r(n)$ define $R(n)$ the number of units adding from the 
right to the inicial configuration $r(0)$ in $n$ steps:
$$
R(n):=\sum\limits_{k=1}^n r_k(k).
$$
We are interested in the behaviour of the process $R(n)$ when $n\rightarrow\infty$.
In the next section we prove low of large numbers, central limit theorem,
local limit theorem, we establish also the large and moderate deviation principles.

Further we suppose that the following condition $[\mathbf{A}]$ holds true. The 
condition $[\mathbf{A}]$ composed by two items.
\emph{\begin{enumerate}
\item Inicial configuration $r(0)$ contains at least one unit.
\item There exist constants $1>\delta_1>\delta_2>0$ such that for all 
$ k\in 0\cup\mathbb{N},$ $1\leq j\leq 2^{v-1}$ the following inequalities holds
$$
\delta_1\geq\mathbf{p}_{kj}\geq\delta_2.
$$
\end{enumerate}}

The condition 1 is an obvious condition for existence of process and for implementation
of transition probabilities. Note, however, that this condition can be omitted adding
the probability $\mathbf{p}_\infty \in (0,1)$ to atribute a character $1$, when 
the sequence consists of only zeros. Condition 2 gives us the possibility of construction 
of arithmetic generalized renewal process which satisfy Cramer moment condition  
$[\mathbf{C}_0]$ and the condition of arithmeticity of $[\mathbf{Z}]$ (see Section 2).

The paper is organized as follows: in Section 2 we introduce our definitions and notations,
we provide the main result Theorem~\ref{t1.5} (low of large numbers, central limit theorem,
local moderate and large limit theorem and principle of moderate large deviations for $R(n)$);
in Section 3 we prove Theorem~\ref{t1.5}; in Section 4 auxiliary lemmas are proved.

\section{Main results, definitions, notations}

To formulate and prove the main result we need some auxiliary processes which 
we define in this section.

For any state $r(n)$ we correspond the pair
$$
Y(n):=(Y_1(n),\mathbf{Y}_2(n)),
$$
where $Y_1(n):=n-m_n$ (the distance to the nearest unit, or, what's the same, the 
number of zeros before the nearest unit),
$$\mathbf{Y}_2(n):=
(r_{m_n-v+1}(n),r_{m_n-v+2}(n), \dots, r_{m_n}(n))
$$
(sequence of the nearest unit and $v-1$ letters from its left).

Note, that the pair $Y(n):=(Y_1(n),\mathbf{Y}_2(n))$ can transit 
with probability 
$\mathbf{p}_{k_{n}j_n}$ into the pair
$$
Y(n+1)=(0,\mathbf{Y}_2(n+1)),
$$
where
$$
\mathbf{Y}_2(n+1)=\left\{
           \begin{array}{ll}
                              (0,\dots,0,1), &

                            \mbox{ if } Y_1(n)\geq v-1,\\
(r_{m_n-v+Y_1(n)+2}(n), \dots, r_{m_n}(n),0,\dots,0,1),& \mbox{ if } Y_1(n)< v-1,\\
                              \end{array}
                              \right.
$$
and with probability $1-\mathbf{p}_{k_{n}j_n}$ into the pair
$$
Y(n+1)=(Y_1(n)+1,\mathbf{Y}_2(n)).
$$

In this way $Y(n+1)$ is an random function on $Y(n)$.
Thus, the sequence $\{Y(n)\}$, $n\geq 0$ is homogeneous Markov chain with phase state
$$
\mathbf{\mathcal{Y}}:=\{y=(y_1,\mathbf{y}_2):y_1\in 0\cup \mathbb{N},\mathbf{y}_2=(a_1,\dots,a_v),
a_1\in \mathcal{A},\dots,a_{v-1}\in \mathcal{A},a_v=1\}.
$$

Let us pick out the state
$$
y_0:=(0,(0,\dots,0,1)).
$$
Note that the chain can jump in one step from any state $(y_1,\mathbf{y}_2)$
to chosen state $y_0$, if the coordinate $y_1$ not less than $v-1$.

Denote
$$
\begin{aligned}
\tau_1&:=\min\{n> 0: Y(n)=y_0\},
\\
\tau_k&:=\min\{n> \tau_1+\dots+\tau_{k-1}: Y(n)=y_0\}-(\tau_1+\dots+\tau_{k-1}), \ \ \ 2\leq k.
\end{aligned}
$$

Since $\{Y(n)\}$ is homogeneous Markov chain, the random variables 
$\tau_1,\dots,\tau_{k},\dots$ are independent and, moreover, 
$\tau_{k}$ are identically distributed when $k\geq 2$.

Let $\zeta_k$ be the number of units added from the right during the time $n \in \{(T_{k-1}+1, \dots,T_k\}$, 
where 
$$
T_0:=0, \ \ \ T_{k}:=T_{k-1}+\tau_k \ \ \ \text{для} \ \ \ k\in \mathbb{N}.
$$
In other words
$$
\zeta_0:=0, \ \ \ \zeta_{k}:=\sum\limits_{n=T_{k-1}+1}^{T_{k}} r_n(n) \ \ \ \text{for} \ \ \ k\in \mathbb{N}.
$$

By construction the random vectors $\xi_k:=(\tau_k,\zeta_k)$, $k\in \mathbb{N}$
are independent, and 
$\xi_k$ are identically distributed when $k\geq 2$.

Let 
$$
\nu(0):=0, \ \ \ \nu(n):=\max\{k\geq 0:T_k<n\} \ \ \ \text{for} \ \ \ k\in \mathbb{N}.
$$

Define \emph{generalized renewal process}
$$
Z(n):=\sum\limits_{k=0}^{\nu(n)}\zeta_k.
$$

Let the random vector $\xi:=(\tau,\zeta)$ has the distribution which coincides 
with distribution of vectors $\xi_k=(\tau_k,\zeta_k)$ for $k\geq 2$.

Since $\zeta_k\leq \tau_k$ a.s. when $k\in \mathbb{N}$, then from Lemma~\ref{l1.1} (see Section 4)
it follows that for $\xi_1$ and $\xi$ the Cramer's condition $[\mathbf{C}_0]$ holds true:
\begin{equation*}
\mathbf{E}e^{\delta|\xi_1|}\leq\mathbf{E}e^{2\delta\tau_1}<\infty \text{ and } 
\mathbf{E}e^{\delta|\xi|}\leq\mathbf{E}e^{2\delta\tau}<\infty, \text{ when } \delta<\rho/2,
\end{equation*}
where $\rho>0$ is the constant from Lemma~\ref{l1.1}.

From Lemma~\ref{l1.1*} (see Section 4) we obtein that the vector $\xi$ satisfies the
arithmeticity condition $[\mathbf{Z}]$:

\begin{enumerate}
 \item[] \emph{For any $u\in \mathbb{Z}^2$ the equality $f(2\pi u)=1$ holds and for any
$u\in \mathbb{R}^2\setminus\mathbb{Z}^2$ the inequality $|f(2\pi u)|<1$ holds,
where for $u=(u_1,u_2)\in\mathbb{R}^2$ the function
$$f(u):=\mathbf{E}e^{i(u_1\tau+u_2\zeta)}$$
is the characteristic function for $\xi$.}
\end{enumerate}
Give the notation we need from the paper \cite{Mog1}

$$
\begin{aligned}
a&:=\frac{\mathbf{E}\zeta}{\mathbf{E}\tau}, \ \ \ \sigma^2:=\frac{1}{a_\tau}\mathbf{E}(\zeta - a \tau)^2,
\text{ где } a_\tau := \mathbf{E}\tau,
\\
\psi_1(\lambda,\mu)&:=\mathbf{E}e^{\lambda\tau_1+\mu\zeta_1}, \ \ \ \psi(\lambda,\mu):=\mathbf{E}e^{\lambda\tau+\mu\zeta},
\\
A(\lambda,\mu)&:=\ln \psi(\lambda,\mu), \ \ \ (\lambda,\mu)\in \mathbb{R}^2,
\\
D(\theta,\alpha)&:=\sup\limits_{(\lambda,\mu)\in \mathcal{A}^{\leq 0}} \{\lambda\theta + \mu \alpha\}, \text{ где }
\mathcal{A}^{\leq 0}:=\{ (\lambda,\mu):A(\lambda,\mu)\leq 0\},
\\
D(\alpha)&:=D(1,\alpha).
\end{aligned}
$$

Denote $\mathfrak{B}$ the Borel $\sigma$-algebra of subsets of $\mathbb{R}$.  
For an arbitrary set $B\in\mathfrak{B}$ we denote 
$[B]$ and $(B)$ its closure and interior correspondingly.

Now we give the definition of the Large Deviation Principle (LDP).
\begin{definition} \label{d1.1} The sequence of random variables $s_n$ satisfies LDP
in  $\mathbb{R}$ with rate function 
$I = I(y)\,:\, \mathbb{R} \rightarrow [0,\infty]$  and normalized function 
$\varphi(n)\,:\, \lim\limits_{n\rightarrow\infty}\varphi(n) = \infty$,
if for any $c \geq 0 $ the set  $\{ y \in \mathbb{R}\,:\, I(y) \leq c \}$ is compact 
and for any set $B \in \mathfrak{B}$ the following inequalities holds true
$$
\begin{aligned}
&\limsup_{n \rightarrow \infty} \frac{1}{\varphi(n)} \ln \mathbf{P}(\, s_n \in B \,) \leq - I([B]), \\
&\liminf_{n \rightarrow \infty} \frac{1}{\varphi(n)} \ln \mathbf{P}(\, s_n \in B \,)\geq -I((B)),
\end{aligned}
$$
where $I(B) = \inf\limits_{y \in B} I(y)$, $I(\varnothing)=\infty$.
\end{definition}

Denote $\Phi_{0,\sigma^2}$ the normal distribution with parameters $(0,\sigma^2)$,
and by $\Rightarrow$ we deonte the convergence in distribution.

Let us give now the main result of our work.

\begin{theorem} \label{t1.5} Let the condition $[\mathbf{A}]$ holds true. Then
\begin{enumerate}
 \item (\textbf{strong low of large numbers}) When $n\rightarrow \infty$
$$
\frac{R(n)}{n}\rightarrow a \text{ a.s.}
$$
\item (\textbf{central limit theorem}) When $n\rightarrow \infty$
$$
\frac{R(n)-an}{\sqrt{n}}\Rightarrow \Phi_{0,\sigma^2}.
$$
\item  (\textbf{local theorem in regions of normal, moderate and large deviations}) 
There exists $\Delta>0$ such that if $x\in 0\cup \mathbb{N}$,
$\lim\limits_{n\rightarrow\infty}\frac{x}{n}=\alpha_0$ and
$|\alpha_0-a|\leq\Delta$, then
$$
\mathbf{P}(R(n)=x)=\frac{1}{\sqrt{n}}\psi_1(\lambda(\alpha_0),\mu(\alpha_0))
C_H(1,\alpha_0)I(\alpha_0)e^{-nD(\frac{x}{n})}(1+o(1)),
$$
where
$$
I(\alpha_0)=\sum\limits_{l=1}^{\infty}e^{\lambda(\alpha_0) l}\mathbf{E}(e^{\mu(\alpha_0) R(l)},\tau\geq l
\mid Y(0)=y_0),
$$
and $C_H(\theta,\alpha)$ is the positive function, which is continuous in a neighborhood of 
the point $(\theta,\alpha)=(1,\alpha_0)$ and it is known explicitly from Theorem 2.1 and 2.1A \cite{Mog1}.

\item  (\textbf{local theorem in regions of normal and moderate deviations}) 
If $x\in 0\cup \mathbb{N}$ and $\lim\limits_{n\rightarrow\infty}\frac{x}{n}=a$,
then the following equality holds true
$$
\mathbf{P}(R(n)=x)=\frac{1}{\sigma\sqrt{2\pi n}}
e^{-nD(\frac{x}{n})}(1+o(1)).
$$
\item (\textbf{moderate deviation principle}) Let  the sequence
$\kappa:=\kappa_n$,  $\kappa\in\mathbb{R}$ satistfies the conditions
$$
\lim\limits_{n\rightarrow\infty}\frac{\kappa}{n}=0, \ \ \ \lim\limits_{n\rightarrow\infty}\frac{\kappa}{\sqrt{n}}=\infty.
$$
Then the sequence of random variables $\tilde{R}(n):=\frac{R(n)-an}{\kappa}$ satisfies LDP
with normalized function $\varphi(n):=\frac{\kappa^2}{n}$ and rate function
$$
I(y):=\frac{y^2}{2\sigma^2}.
$$
\end{enumerate}
\end{theorem}

\section{Proof of Theorem~\ref{t1.5}}
\DD \ of statements 1) and 2). Since $\zeta_k\leq \tau_k$ a.s. when $k\in \mathbb{N}$, then
the following inequality holds
\beq \label{2}
Z(n)\leq R(n)\leq Z(n)+\tau_{\nu(n)+1}.
\eeq
Using Lemma~\ref{l1.21} and Borel-Cantelli lemma, when $n\rightarrow\infty$, it is easy to see that 
$$
\frac{\tau_{\nu(n)+1}}{\sqrt{n}}\rightarrow 0 \text{ a.s.}
$$
Thus the statements 1) and 2) follow from the inequality (\ref{2}) and corresponding 
results for $Z(n)$ 
(see \cite{Bor1} Theorem 11.5.2 pp.332 and Theorem 10.6.2 pp. 311).
\\$\Box$

\bigskip

\DD \ of statement 3). Consider
\beq \label{1.12}
\begin{aligned}
L_n(x) &:=\sum\limits_{k=1}^\infty\sum\limits_{l=1}^{[\ln^2n]}\mathbf{P}(R(n)=x,T_k=n-l,\tau_{k+1}\geq l)
\\
&=\sum\limits_{k=1}^\infty\sum\limits_{l=1}^{[\ln^2n]}\sum\limits_{s=0}^{l}
\mathbf{P}(Z_k=x-s,R(n)-Z_k=s,T_k=n-l,\tau_{k+1}\geq l).
\end{aligned}
\eeq
Note that if $T_k=n-l$, then the random variable $Z_k$ uniquely defined by the values 
of Markov chain $Y(m)$ when $m<n-l$, but the random variables
$R(n)-Z_k$ and $\tau_{k+1}$ depend on the values of the chain $Y(m)$ when $m>n-l$. 
Therefore, by the inclusion 
$$
\{\omega:T_k=n-l\}\subseteq\{\omega:Y(n-l)=y_0\}
$$
the following equality holds
\beq \label{1.13}
\begin{aligned}
\mathbf{P} & (Z_k=x-s,R(n)-Z_k=s,\tau_{k+1}\geq l \mid T_k=n-l)
\\
&=\mathbf{P}(Z_k=x-s \mid T_k=n-l)\mathbf{P}(R(n)-Z_k=s,\tau_{k+1}\geq l \mid T_k=n-l).
\end{aligned}
\eeq
Applying (\ref{1.12}) and (\ref{1.13}) we obtain
\beq \label{1.14}
\begin{aligned}
& L_n(x)=\sum\limits_{k=1}^\infty\sum\limits_{l=1}^{[\ln^2n]}\sum\limits_{s=0}^{l}
\mathbf{P}(Z_k=x-s,R(n)-Z_k=s,T_k=n-l,\tau_{k+1}\geq l)
\\
&=\sum\limits_{k=1}^\infty\sum\limits_{l=1}^{[\ln^2n]}\sum\limits_{s=0}^{l}
\mathbf{P}(Z_k=x-s|T_k=n-l)\mathbf{P}(R(n)-Z_k=s,\tau_{k+1}\geq l|T_k=n-l)\mathbf{P}(T_k=n-l)
\\
&=\sum\limits_{k=1}^\infty\sum\limits_{l=1}^{[\ln^2n]}\sum\limits_{s=0}^{l}
\mathbf{P}(Z_k=x-s,T_k=n-l)\mathbf{P}(R(n)-Z_k=s,\tau_{k+1}\geq l|T_k=n-l).
\end{aligned}
\eeq
Note that
$$
\mathbf{P}(R(n)-Z_k=s,\tau_{k+1}\geq l \mid T_k=n-l)=\mathbf{P}(R(l)=s,\tau\geq l \mid Y(0)=y_0).
$$
Thus, from equality (\ref{1.14}) it is follows that
\beq \label{1.15}
L_n(x)=\sum\limits_{l=1}^{[\ln^2n]}\sum\limits_{s=0}^{l}\mathbf{P}(R(l)=s,\tau\geq l|Y(0)=y_0)
\sum\limits_{k=1}^\infty\mathbf{P}(Z_k=x-s,T_k=n-l).
\eeq
Since $\mathbf{P}(T_0=n-l)=0$ when $0  \leq l \leq [\ln^2n]$, then from Theorem~2.2А \cite{Mog1}
it is follows, when $n\rightarrow\infty$, that
$$
\begin{aligned}
\sum\limits_{k=1}^\infty  \mathbf{P}(Z_k=x-s,T_k=n-l)=\sum\limits_{k=0}^\infty\mathbf{P}(Z_k=x-s,T_k=n-l)&
\\
=\frac{1}{\sqrt{n}}\psi_1\bigg(\lambda\bigg(\frac{\tilde{\alpha}}{\theta}\bigg),\mu\bigg(\frac{\tilde{\alpha}}{\theta}\bigg)\bigg)
C_H(\theta,\tilde{\alpha})e^{-nD(\theta,\tilde{\alpha})}(1+o(1))&,
\end{aligned}
$$
where
$$
\tilde{\alpha}:=\frac{x-s}{n}, \ \ \ \theta:=\frac{n-l}{n}.
$$
Since the function $\psi_1(\lambda(\alpha),\mu(\alpha))$ is continuous in a neighborhood of the point
$\alpha=a$, and the function $C_H(\theta,\alpha)$ is continuous in a neighborhood of the point
$(\theta,\alpha)=(1,a)$, then, for sufficiently small $\Delta$ and $n\rightarrow\infty$ 
the following equality holds true
\beq \label{1.151}
\sum\limits_{k=1}^\infty\mathbf{P}(Z_k=x-s,T_k=n-l)
=\frac{1}{\sqrt{n}}\psi_1(\lambda(\alpha_0),\mu(\alpha_0))
C_H(1,\alpha_0)e^{-nD(\theta,\tilde{\alpha})}(1+o(1)).
\eeq
Applying Lemma~\ref{l1.4} (see Section 4) and considering that $0 \leq s \leq l \leq [\ln^2n]$
and $|\alpha_0-a|<\Delta$, from the equality (\ref{1.151}) we obtain
\beq \label{1.16}
\begin{aligned}
&\sum\limits_{k=1}^\infty\mathbf{P}(Z_k=x-s,T_k=n-l)
\\
&=\frac{1}{\sqrt{n}}\psi_1(\lambda(\alpha_0),\mu(\alpha_0))
C_H(1,\alpha_0)e^{-nD(\frac{x}{n})+(\lambda(\frac{x}{n})+\varepsilon_n) l
+(\mu(\frac{x}{n})+\theta_n) s}(1+o(1)).
\end{aligned}
\eeq
Show that 
\beq \label{1.161}
\begin{aligned}
\lim\limits_{n\rightarrow\infty} & \sum\limits_{l=1}^{[\ln^2n]}\sum\limits_{s=0}^{l}\mathbf{P}(R(l)=s,\tau\geq l|Y(0)=y_0)
e^{(\lambda(\frac{x}{n})+\varepsilon_n) l
+(\mu(\frac{x}{n})+\theta_n) s}
\\
&=\sum\limits_{l=1}^{\infty}e^{\lambda(\alpha_0) l}\mathbf{E}(e^{\mu(\alpha_0) R(l)},\tau\geq l|Y(0)=y_0)=:I(\alpha_0).
\end{aligned}
\eeq
Due to the fact that
$$
\lim\limits_{n\rightarrow\infty}e^{(\lambda(\frac{x}{n})+\varepsilon_n) l
+(\mu(\frac{x}{n})+\theta_n) s}=e^{\lambda(\alpha_0) l
+\mu(\alpha_0) s}$$
the equality (\ref{1.161}) will be proved if we can show that the series 
\beq \label{1.162}
\sum\limits_{l=1}^{\infty}\sum\limits_{s=0}^{l}\mathbf{P}(R(l)=s,\tau\geq l|Y(0)=y_0)
e^{(\lambda(\frac{x}{n})+\varepsilon_n) l
+(\mu(\frac{x}{n})+\theta_n) s}
\eeq
converges.

Note that if $\tau\geq l$, then $\zeta\geq R(l)$, thus 
\beq\label{1.17}
\begin{aligned}
\sum\limits_{l=1}^{\infty}\sum\limits_{s=0}^{l} & \mathbf{P}(R(l)=s,\tau\geq l \mid Y(0)=y_0)
e^{(\lambda(\frac{x}{n})+\varepsilon_n) l
+(\mu(\frac{x}{n})+\theta_n) s}
\\
&=\sum\limits_{l=1}^{\infty}e^{(\lambda(\frac{x}{n})+\varepsilon_n)l}
\mathbf{E}(e^{(\mu(\frac{x}{n})+\theta_n)R(l)},\tau\geq l\mid Y(0)=y_0)
\\
&\leq\sum\limits_{l=1}^{\infty}e^{(\lambda(\frac{x}{n})+\varepsilon_n) l}
\mathbf{E}(e^{(\mu(\frac{x}{n})+\theta_n)\zeta},\tau\geq l \mid Y(0)=y_0).
\end{aligned}
\eeq
Due to the Cramer's condition $[\mathbf{C}_0]$ for sufficiently small $\Delta>0$ and sufficiently
large $n$
$$
\mathbf{E}(e^{2(\mu(\frac{x}{n})+\theta_n)\zeta} \mid Y(0)=y_0)<\infty
$$
and there exists $\rho>0$ such that
$$
\mathbf{E}(e^{2(\lambda(\frac{x}{n})+\varepsilon_n+\rho)\tau} \mid Y(0)=y_0)<\infty.
$$
Therefore, using Cauchy-Bunyakovsky and Chebyshev inequalities, we obtain
\beq \label{1.171}
\begin{aligned}
\mathbf{E} & (e^{(\mu(\frac{x}{n})+\theta_n)\zeta},\tau\geq l \mid Y(0)=y_0)
\\
& \leq
(\mathbf{E}e^{2(\mu(\frac{x}{n})+\theta_n)\zeta} \mid Y(0)=y_0)^\frac{1}{2}
(\mathbf{P}(\tau\geq l \mid Y(0)=y_0))^\frac{1}{2}
\\
&\leq (\mathbf{E}e^{2(\mu(\frac{x}{n})+\theta_n)\zeta} \mid Y(0)=y_0)^\frac{1}{2}
(\mathbf{E}e^{2(\lambda(\frac{x}{n})+\varepsilon_n+\rho)\tau} \mid Y(0)=y_0)^\frac{1}{2}
e^{-(\lambda(\frac{x}{n})+\varepsilon_n+\rho)l}
\\
&
=:Ke^{-(\lambda(\frac{x}{n})+\varepsilon_n+\rho)l}.
\end{aligned}
\eeq
Using (\ref{1.17}), (\ref{1.171}), we obtain
$$
\sum\limits_{l=1}^{\infty}\sum\limits_{s=0}^{l}\mathbf{P}(R(l)=s,\tau\geq l \mid Y(0)=y_0)
e^{(\lambda(\frac{x}{n})+\varepsilon_n) l
+(\mu(\frac{x}{n})+\theta_n) s}\leq
K\sum\limits_{l=1}^{\infty}e^{-\rho l}<\infty.
$$
Thus, the series (\ref{1.162}) converges, hence the equality (\ref{1.161}) holds true.

From (\ref{1.15}), (\ref{1.16}) and (\ref{1.161}) it is follows that
\beq \label{7.1}
L_n(x)=\frac{1}{\sqrt{n}}\psi_1(\lambda(\alpha_0),\mu(\alpha_0))
C_H(1,\alpha_0)I(\alpha_0)e^{-nD(\frac{x}{n})}(1+o(1)).
\eeq
It is obvious that the following inequality holds
\beq \label{7.2}
L_n(x)\leq \mathbf{P}(R(n)=x) \leq L_n(x)+\mathbf{P}(R(n)=x,\tau_{\nu(n)+1}\geq \ln^2n).
\eeq
From Lemma~\ref{l1.3} (see Section 4) it is follows that
\beq \label{7.3}
\mathbf{P}(R(n)=x,\tau_{\nu(n)+1}\geq \ln^2n)\leq\tilde{C} e^{- nD(\alpha)-\tilde{\gamma}\ln^2n}.
\eeq

From (\ref{7.1}), (\ref{7.2}) and (\ref{7.3}) it is follows that
\beq\label{1.1711}
\begin{aligned}
\mathbf{P}(R(n)=x) & = (1+o(1))L_n(x)
\\
& =\frac{1}{\sqrt{n}}\psi_1(\lambda(\alpha_0),\mu(\alpha_0))
C_H(1,\alpha_0)I(\alpha_0)e^{-nD(\frac{x}{n})}(1+o(1)).
\end{aligned}
\eeq
$\Box$

\bigskip

\DD \ of statement 4). Due to the fact that $(\lambda(a),\mu(a))=(0,0)$ and the function $I(\alpha)$
is continuous in a neighborhood of the point $\alpha=a$ we will have 
\beq \label{1.18}
\begin{aligned}
I(a) &=\sum\limits_{l=1}^{\infty}e^{\lambda(a) l}\mathbf{E}(e^{\mu(a) R(l)},\tau\geq l \mid Y(0)=y_0)
\\
&=\sum\limits_{l=1}^{\infty}\mathbf{P}(\tau\geq l \mid Y(0)=y_0)=\mathbf{E}(\tau \mid Y(0)=y_0).
\end{aligned}
\eeq
From Lemma~2.1 (see \cite{Mog1}) it is follows that
$$
C_H(1,a)=\frac{1}{\mathbf{E}(\tau \mid Y(0)=y_0)\sigma\sqrt{2\pi}}.
$$
Hence, from (\ref{1.1711}) and (\ref{1.18}) we obtain
$$
\mathbf{P}(R(n)=x)=\frac{1}{\sigma\sqrt{2\pi n}}
e^{-nD(\frac{x}{n})}(1+o(1)).
$$
$\Box$

\DD \ of statement  5). From Consequence~3.2 (see \cite{LogMog}) it is follows that 
the sequence of random variables $\tilde{Z}(n):=\frac{Z(n)-an}{\kappa}$ satisfies
LDP with normalized function $\varphi(n)=\frac{\kappa^2}{n}$ and rate funtion $I(y)$.

Using Lemma~\ref{l1.21}, for any $\varepsilon>0$ we will have
$$
\lim_{n\rightarrow\infty}\frac{n}{\kappa^2}\ln \mathbf{P}(|\tilde{R}(n)-\tilde{Z}(n)|>\varepsilon)\leq
\lim_{n\rightarrow\infty}\frac{n}{\kappa^2}\ln \mathbf{P}(\tau_{\nu(n)+1}>\kappa\varepsilon)
$$
$$
\leq \lim_{n\rightarrow\infty}\frac{n}{\kappa^2}\ln e^{-\frac{\rho}{4}\kappa\varepsilon}=
-\lim_{n\rightarrow\infty}\frac{n\rho\varepsilon}{4\kappa}=-\infty.
$$
Therefore from Theorem 4.2.13 (see \cite{Demb}) we obtain that the sequences $\tilde{R}(n)$ and
$\tilde{Z}(n)$ satisfy the same LDP. $\Box$

\section{Auxiliary Results}

\begin{lemma} \label{l1.1} For any 
$k, \ n\in \mathbb{N}$ the following inequality holds
\beq \label{1.2}
\mathbf{P}(\tau_{k}\geq n)\leq C e^{-\rho n},
\eeq
where
$$
C:=\frac{1}{1-(1-\delta_1)^{v-1}\delta_2}, \ \ \
\rho:=\frac{1}{v}\ln C.
$$
\end{lemma}

\DD. Due to the fact that the process $Y(n)$ is markovian it is sufficient to prove 
Lemma~\ref{l1.1} for $\tau_{1}$ with an arbitrary initial condition.
We fix some inicial state $Y(0)=(y_1,\mathbf{y}_2)$.
Since $C>1$, $Ce^{-\rho n}=C^{1-\frac{n}{v}}$, then for $n\leq v$ the right-hand side of
inequality (\ref{1.2}) is not less than 1, therefore (\ref{1.2}) obviously holds true.

We prove now the inequality (\ref{1.2}) when $n\geq v+1$. Denote
$k:=[\frac{n}{v}]$ and for $l=0, 1, \dots, k-1$
we consider the events
$$
A_l:=\{\omega:r_{vl+1}(n)=0,\dots,r_{vl+v-1}(n)=0,r_{vl+v}(n)=1\}.
$$

Denote $B$ the complement to the set $\bigcup\limits_{l=0}^{k-1}A_l$:
$$
B:=\bigcap\limits_{l=0}^{k-1}\overline{A_l}.
$$
Since it is obvious that
$$
\{\tau_1<n\}\supset \bigcup\limits_{l=0}^{k-1}A_l,
$$
then we obtain 
$$
\{\tau_1\geq n\}\subseteq B.
$$
Hence we have
$$
\begin{aligned}
P_n:&=\mathbf{P}(\tau_1\geq n \  \mid  \ Y(0)=(y_1,\mathbf{y}_2))
\\
&\leq \mathbf{P}(B \  \mid  \ Y(0)=(y_1,\mathbf{y}_2))
= \mathbf{P}\bigg(\bigcap\limits_{l=0}^{k-1}\overline{A_l} \ \bigg| \ Y(0)=(y_1,\mathbf{y}_2)\bigg)
\\
&=\mathbf{P}\left(\overline{A_0} \ \mid \ Y(0)=(y_1,\mathbf{y}_2)\right)
\prod\limits_{l=1}^{k-1}\mathbf{P}\bigg(\overline{A_l} \ \bigg|
\bigcap\limits_{i=0}^{l-1}\overline{A_i}, \ Y(0)=(y_1,\mathbf{y}_2)\bigg).
\end{aligned}
$$
Since by condition $[\mathbf{A}]$ each cofactor in the right-hand side has 
the following upper bound
$$
1-(1-\delta_1)^{v-1}\delta_2=\frac{1}{C},
$$
then we have
$$
P_n\leq \left(\frac{1}{C}\right)^{k}
\leq \left(\frac{1}{C}\right)^{\frac{n}{v}-1}=Ce^{-\rho n}.
$$
$\Box$

\begin{lemma} \label{l1.1*} For any $u\in \mathbb{Z}^2$ the equality $f(2\pi u)=1$ holds true, and
for any $u\in \mathbb{R}^2\setminus\mathbb{Z}^2$ the inequality $|f(2\pi u)|<1$ holds, where
$$f(u):=\mathbf{E}e^{i(u_1\tau+u_2\zeta)}$$
is characteristic function for $\xi$.
\end{lemma}

\DD. Since $\tau$ and $\zeta$ are integers numbers, then, it is obvious, that 
for $u\in \mathbb{Z}^2$ the equality 
$f(2\pi u)=1$ holds true. We show that for any 
$u\in \mathbb{R}^2\setminus\mathbb{Z}^2$ the inequality $|f(2\pi u)|<1$ holds true.

Suppose that it is not true, then there exists $(u_1,u_2)\in\mathbb{R}^2\setminus\mathbb{Z}^2$
such that $|\mathbf{E}e^{2\pi i(u_1\tau+u_2\zeta)}|=1$.
Note that equality $|\mathbf{E}e^{2\pi i(u_1\tau+u_2\zeta)}|=1$ is equivalent to the fact that
there exists $k\in \mathbb{R}$ such that
$$2\pi(u_1\tau+u_2\zeta)=k  \mod(2\pi) \text{ a.s.}$$
From Condition $[\mathbf{A}]$ it is follows that
$$
\begin{aligned}
&\mathbf{P}(\zeta=1,\tau=s+1)\geq (1-\delta_1)^s\delta_2>0,\\
&\mathbf{P}(\zeta=2,\tau=s+1)\geq \delta_2(1-\delta_1)^{s-1}\delta_2>0,\\
&\mathbf{P}(\zeta=1,\tau=s+2)\geq (1-\delta_1)^{s+1}\delta_2>0.
\end{aligned}
$$
Thus, if or hypothesis is true, then should exist $k_1\in \mathbb{Z}$, $k_2\in \mathbb{Z}$,
$k_3\in \mathbb{Z}$ such that the following inequalities holds true
$$
\left\{ \begin{array}{ll} 2\pi(u_1(s+1)+u_2)=k+2\pi k_1\\
2\pi(u_1(s+1)+2u_2)=k+2\pi k_2\\
2\pi(u_1(s+2)+u_2)=k+2\pi k_3.
\end{array} \right. \label{8a}
$$
Divide each equality by $2\pi$. Subtracting from the 2nd equality the 1st we obtain
$u_2=k_2-k_1\in \mathbb{Z}$; substracting from the 3rd equality the 1st, we obtain
$u_1=k_3-k_1\in \mathbb{Z}$. The resulting contradiction
completes the proof.
$\Box$

\bigskip

For the vector $(\tilde{\lambda},\tilde{\mu})$ such that $\psi(\tilde{\lambda},\tilde{\mu})=1$, 
we consider the sequence of random vectors
$(\hat{\tau}_k,\hat{\zeta}_k)$, $k\in\mathbb{N}$, whose joint distribution is given as follows
\beq\label{1.5}
\mathbf{P}((\hat{\tau}_1,\hat{\zeta}_1)\in A_1,\dots,(\hat{\tau}_k,\hat{\zeta}_k)\in A_k,\dots):=
\frac{1}{\psi_1(\tilde{\lambda},\tilde{\mu})}
\prod\limits_{k=1}^\infty\mathbf{E}(e^{\tilde{\lambda}\tau_k+\tilde{\mu}\zeta_k};(\tau_k,\zeta_k)\in A_k).
\eeq
Let $\hat{\tau}_0:=0$, $\hat{\zeta}_0:=0$,  $\hat{\nu}(0):=0$.
Denote
$$
\hat{T}_k:=\sum\limits_{l=0}^k\hat{\tau}_l, \ \ \ \hat{\nu}(n):=\max\{k\geq 0:\hat{T}_k<n\}.
$$

\begin{lemma} \label{l1.2} Let $\gamma+\tilde{\lambda}+\tilde{\mu}<\rho$, then there exists
the constant $\hat{C}>0$,
such that for any $n\in \mathbb{N}$ the following inequality holds true
$$
\mathbf{E} e^{\gamma \hat{\tau}_{\hat{\nu}(n)+1}}<\hat{C}n.
$$
\end{lemma}

\DD. Since random variables $\hat{\tau}_{k+1}$ и $\hat{T}_k$ are independent, then
$$
\begin{aligned}
E_1:=\mathbf{E} e^{\gamma \hat{\tau}_{\hat{\nu}(n)+1}} &=\mathbf{E} (e^{\gamma \hat{\tau}_{1}};\hat{\tau}_1\geq n)+
\sum\limits_{k=1}^\infty\mathbf{E} (e^{\gamma \hat{\tau}_{k+1}};\hat{T}_k< n\leq \hat{T}_{k+1})
\\
&\leq \frac{1}{\psi_1(\tilde{\lambda},\tilde{\mu})}\mathbf{E} e^{\gamma \tau_{1}} + \mathbf{E} e^{\gamma \tau+\tilde{\lambda}\tau+\tilde{\mu}\zeta}
\sum\limits_{k=1}^\infty\mathbf{P} (\hat{T}_k< n).
\end{aligned}
$$
Due to the fact that by arithmeticity the inequality $\hat{T}_k\geq n$ a.s. when $k\geq n$.
Therefore, using Lemma~\ref{l1.1} and inequality $\tau_k\geq \zeta_k$ a.s., we obtain
$$
\begin{aligned}
E_1&\leq \frac{1}{\psi_1(\tilde{\lambda},\tilde{\mu})}\mathbf{E} e^{\gamma \tau_{1}}
 +  \mathbf{E} e^{(\gamma+\tilde{\lambda}+\tilde{\mu})\tau}n
\leq \frac{1}{\psi_1(\tilde{\lambda},\tilde{\mu})}C\sum\limits_{k=1}^\infty  e^{(\gamma-\rho) k} + Cn\sum\limits_{k=1}^\infty  e^{(\gamma+\tilde{\lambda}+\tilde{\mu}-\rho) k}
\\
&\leq \left(\frac{C}{\psi_1(\tilde{\lambda},\tilde{\mu})}\frac{ e^{(\gamma-\rho)}}{1-e^{(\gamma-\rho)}}+
\frac{Ce^{(\gamma+\tilde{\lambda}+\tilde{\mu}-\rho)}}{1-e^{(\gamma+\tilde{\lambda}+\tilde{\mu}-\rho)}}\right)n.
\end{aligned}
$$
$\Box$

\begin{lemma} \label{l1.21} Let
$\lim\limits_{n\rightarrow\infty}\kappa_n=\infty$.
Then for sufficiently large $n\in \mathbb{N}$
the following holds true
$$
\mathbf{P}(\tau_{\nu(n)+1}\geq \kappa_n)\leq e^{-\frac{\rho}{4}\kappa_n}.
$$
\end{lemma}

\DD. Since random variables $\tau_{k+1}$ and $T_k$ are independent, then
$$
\begin{aligned}
\mathbf{P}(\tau_{\nu(n)+1} \geq \kappa_n) &\leq\mathbf{P}(\tau_1\geq \kappa_n)+
\sum\limits_{k=1}^\infty\mathbf{P} (\tau_{k+1}\geq \kappa_n,T_k< n\leq T_{k+1})
\\
&\leq\mathbf{P}(\tau_1\geq \kappa_n)+
\mathbf{P} (\tau\geq \kappa_n)\sum\limits_{k=1}^\infty\mathbf{P} (T_k< n).
\end{aligned}
$$
Due to the arithmeticity the inequality $T_k\geq n$ holds true almost surely when $k\geq n$.
By Lemma~\ref{l1.1} and Chebyshev inequality for sufficiently large $n$ we obtain
$$
\mathbf{P}(\tau_{\nu(n)+1}\geq\kappa_n)\leq\frac{\mathbf{E}e^{\frac{\rho}{2}\tau_1}}{e^{\frac{\rho}{2}\kappa_n}}+
n\frac{\mathbf{E}e^{\frac{\rho}{2}\tau}}{e^{\frac{\rho}{2}\kappa_n}}\leq e^{-\frac{\rho}{4}\kappa_n}.
$$
$\Box$

\begin{lemma} \label{l1.3} There exist constants
$\Delta>0$, $\tilde{C}>0$, $\tilde{\gamma}>0$ such that for 
$x\in 0\cup \mathbb{N}$, $\alpha:=\frac{x}{n}$, $n\geq 1$, $|\alpha-a|\leq\Delta$
the following inequality holds true
$$
\mathbf{P}(R(n)=x,\tau_{\nu(n)+1}\geq \ln^2n)\leq \tilde{C} e^{- nD(\alpha)-\tilde{\gamma}\ln^2n}.
$$
\end{lemma}

\DD. By Theorem~2.1, \cite{Mog1}, it is follows that for sufficiently small $\Delta$ there exists
$\lambda(\alpha)$ and $\mu(\alpha)$ such that $(\lambda(\alpha),\mu(\alpha))\in \mathcal{A}^{\leq 0}$,
$A(\lambda(\alpha),\mu(\alpha))=0$ and
$$
D(\alpha)=\lambda(\alpha)+\alpha\mu(\alpha).
$$
Denote
$$
B_n:=\{\omega:\tau_{\nu(n)+1}\geq \ln^2n\}.
$$
We have
$$
\begin{aligned}
&\mathbf{P}  (R(n)=x,\tau_{\nu(n)+1}\geq \ln^2n)
\\
& =\mathbf{P}(R(n)=x,B_n,\nu(n)=0)
+
\mathbf{P}(R(n)=x,B_n,\nu(n)\geq 1):=\mathbf{P}_0+\mathbf{P}_1.
\end{aligned}
$$
From Lemma~\ref{l1.1} it is follows that
$$
\mathbf{P}_0\leq \mathbf{P}(\tau_1\geq n)\leq C e^{-\rho n}.
$$
Since the function $D(\alpha)$ is continuouis in a neighborhood of the point
$\alpha=a$ and $D(a)=0$, then for sufficiently small $\Delta>0$
and $\alpha: \ |\alpha-a|\leq \Delta$ the following inequality holds
\beq \label{1.7}
\mathbf{P}_0\leq  C e^{-\rho n}\leq C e^{-D(\alpha) n-c\ln^2 n}.
\eeq

Denote $Z_k:=\sum\limits_{l=1}^k\zeta_l$.

Estimate from above $\mathbf{P}_1$. For $\lambda=\lambda(\alpha)$, $\mu=\mu(\alpha)$ we obtain
\beq \label{1.8}
\begin{aligned}
& \mathbf{P}_1=\mathbf{P}(R(n)=x,B_n,\nu(n)\geq 1)
=\sum\limits_{k=1}^\infty\mathbf{P}(R(n)=x,B_n,\nu(n)=k)
\\
&=\sum\limits_{k=1}^\infty\mathbf{P}(R(n)=x,\tau_{k+1}\geq \ln^2 n,T_k<n\leq T_{k+1})
\\
&=e^{-D(\alpha) n}\sum\limits_{k=1}^\infty
\mathbf{E}(e^{-\lambda (T_{k+1}-n)-\mu (Z_{k+1}-x)+\lambda T_{k+1}+\mu Z_{k+1}}; R(n)=x,\tau_{k+1}\geq \ln^2 n,
T_k<n\leq T_{k+1}).
\end{aligned}
\eeq

Note that if $T_{k+1}\geq n$ and $R(n)=x$, then
$x+\zeta_{k+1}\geq Z_{k+1}$, therefore from (\ref{1.8}) it is follows that
\beq \label{1.9}
\begin{aligned}
\mathbf{P}_1 &\leq e^{-D(\alpha) n}\sum\limits_{k=1}^\infty
\mathbf{E}(e^{|\lambda|\tau_{k+1}+|\mu| \zeta_{k+1}+\lambda T_{k+1}+\mu Z_{k+1}};  R(n)=x,\tau_{k+1}\geq \ln^2 n,T_k<n\leq T_{k+1})
\\
&\leq e^{-D(\alpha) n}\sum\limits_{k=1}^\infty
\mathbf{E}(e^{|\lambda| \hat{\tau}_{k+1}+|\mu| \hat{\zeta}_{k+1}};\hat{\tau}_{k+1}\geq \ln^2 n,\hat{T}_k<n\leq \hat{T}_{k+1}),
\end{aligned}
\eeq
where the joint distribution of random variables $(\hat{\tau}_k,\hat{\zeta}_k)$, $k\in\mathbb{N}$ 
has the form
(compare with (\ref{1.5}))
\beq \label{1.10}
\mathbf{P}((\hat{\tau}_1,\hat{\zeta}_1)\in A_1,\dots,(\hat{\tau}_k,\hat{\zeta}_k)\in A_k,\dots):=
\frac{1}{\psi_1(\lambda(\alpha),\mu(\alpha))}
\prod\limits_{k=1}^\infty\mathbf{E}(e^{\lambda(\alpha)\tau_k+\mu(\alpha)\zeta_k};(\tau_k,\zeta_k)\in A_k),
\eeq
Making the summation in the inequality (\ref{1.9}) we obtain
$$
\mathbf{P}_1\leq e^{-D(\alpha) n}
\mathbf{E}(e^{|\lambda| \hat{\tau}_{\nu(n)+1}+|\mu| \hat{\zeta}_{\nu(n)+1}};\hat{\tau}_{\nu(n)+1}\geq \ln^2 n).
$$

Since $\tau_k\geq\zeta_k$ a.s., then from (\ref{1.10}) it is follows that $\hat{\tau}_k\geq\hat{\zeta}_k$ a.s.
Therefore, using Cauchy-Schwarz-Bunjakowski inequality we obtain
\beq \label{1.11}
\mathbf{P}_1\leq e^{-D(\alpha) n}\left(\mathbf{E}e^{2(|\lambda|+|\mu|) \hat{\tau}_{\nu(n)+1}}\right)^\frac{1}{2}
\left(\mathbf{P}(\hat{\tau}_{\nu(n)+1}\geq \ln^2 n)\right)^{\frac{1}{2}}.
\eeq
From Consequence~2.1, \cite{Mog1}, it is follows that $(\lambda(a),\mu(a))=(0,0)$, 
therefore due to the continuity in the point $\alpha=a$ of the function 
$D(\alpha)$ for sufficiently small $\Delta>0$ the following inequality 
$3(|\lambda(\alpha)|+|\mu(\alpha)|)<\rho$ holds true.

Thus, from Lemma~\ref{l1.2} it is follows that there exists a constant $C_1>0$ such that
$$
\mathbf{E}e^{2(|\lambda|+|\mu|) \hat{\tau}_{\nu(n)+1}}\leq C_1 n.
$$

Using Lemma~\ref{l1.2} and Chebyshev inequality for some $C_2>0$, $\gamma>0$ we obtain
$$
\mathbf{P}(\hat{\tau}_{\nu(n)+1}\geq \ln^2 n)<C_2e^{-\gamma \ln^2 n}.
$$
From inequalities (\ref{1.7}), (\ref{1.11}) it is follows that for sufficiently large $n$
$$
\begin{aligned}
\mathbf{P}(R(n)=x,\tau_{\nu(n)+1}\geq \ln^2n) &\leq
C e^{-D(\alpha) n-c\ln^2 n}+e^{-D(\alpha) n}\sqrt{C_1C_2}\sqrt{n}e^{-\frac{\gamma}{2} \ln^2 n}
\\
&\leq C e^{-D(\alpha) n-c\ln^2 n}+\sqrt{C_1C_2}e^{-D(\alpha) n-\frac{\gamma}{4} \ln^2 n}.
\end{aligned}
$$
$\Box$

\begin{lemma} \label{l1.4} There exists $\Delta>0$ such that if
$\lim\limits_{n\rightarrow\infty}\frac{x}{n} = \alpha_0$, $|a-\alpha_0|\leq\Delta$,
then the following holds
\beq \label{1.111}
-nD\bigg(1-\frac{m}{n},\frac{x}{n}-\frac{y}{n}\bigg)=-nD\bigg(1,\frac{x}{n}\bigg)+
\bigg(\lambda\bigg(\frac{x}{n}\bigg)+\varepsilon_n\bigg) m
+\bigg(\mu\bigg(\frac{x}{n}\bigg)+\theta_n\bigg) y,
\eeq
where functions $\varepsilon_n=\varepsilon_n(m,y)$, $\theta_n=\theta_n(m,y)$
satisfy 
$$
\beta_n:=\max\limits_{(m,y)\in \mathcal{B}_n}\{|\varepsilon_n(m,y)|+|\theta_n(m,y)|\}=o(1) \ \ \
\text{when} \ \ \ n\rightarrow \infty,
$$
where $ \mathcal{B}_n:=\{(m,y)\in \mathbb{Z}^2: 1\leq m \leq [n\kappa_n], 1\leq |y| \leq [n\kappa_n]\}$,
 $\kappa_n=o(1)$ when $n\rightarrow \infty$.
\end{lemma}

\DD. For sufficiently small $\Delta$ if $|a-\alpha_0|\leq\Delta$, then the function $D(\alpha)$ 
is analytic in some neighborhood of the point
$\alpha_0$. Therefore, due to the fact that $D(\theta,\alpha)=\theta D(\frac{\alpha}{\theta})$
the function $D(\theta,\alpha)$ is analytic in some neighborhood of the point $(1,\alpha_0)$.
It means that for sufficiently large $n$ the function can be represented as a Taylor series
in a neighborhood of the point $(1,\frac{x}{n})$. Therefore we have
\beq \label{1.112}
-nD\bigg(1-\frac{m}{n},\frac{x}{n}-\frac{y}{n}\bigg)=-nD\bigg(1,\frac{x}{n}\bigg)+\lambda\bigg(\frac{x}{n}\bigg) m
+\mu\bigg(\frac{x}{n}\bigg)y-nM_1,
\eeq
where $M_1$ is the remainder term in Lagrange's form

Denote
$$
D''_{(1,1)}(x,y):=\frac{\partial^2}{\partial^2x}D(x,y), \ \ \ D''_{(1,2)}(x,y):=\frac{\partial^2}{\partial x\partial y}D(x,y),
 \ \ \ D''_{(2,2)}(x,y):=\frac{\partial^2}{\partial^2y}D(x,y).
$$
Then there exist $u\in[1-\frac{m}{n},1]$, $u\in[\frac{x}{n}-\frac{y}{n},\frac{x}{n}]$ such that
$$
|M_1|\leq 2\max(|D''_{(1,1)}(u,v)|,|D''_{(1,2)}(u,v)|,|D''_{(1,2)}(u,v)|)
\bigg(\frac{m^2}{n^2}+\frac{y^2}{n^2}\bigg):=K\bigg(\frac{m^2}{n^2}+\frac{y^2}{n^2}\bigg).
$$
Therefore, from (\ref{1.112}) it is follows that there exist $\varepsilon_n(m,y)$,
$\theta_n(m,y)$ such that
$$
|\varepsilon_n(m,y)|\leq K \frac{m}{n}, \ \ \ |\theta_n(m,y)|\leq K \frac{y}{n}
$$
and the equality (\ref{1.111}) holds.
\\$\Box$

\end{document}